%&amstex          
\input amstex\documentstyle{amsppt}  
\pagewidth{12.5cm}\pageheight{19cm}\magnification\magstep1
\topmatter
\title On the trace of Coxeter elements\endtitle
\author G. Lusztig\endauthor
\address{Department of Mathematics, M.I.T., Cambridge, MA 02139}\endaddress
\thanks{Supported by NSF grant DMS-2153741}\endthanks
\endtopmatter   
\document
\define\mat{\pmatrix}
\define\endmat{\endpmatrix}

\define\Irr{\text{\rm Irr}}

\define\si{\sim}

\define\sqc{\sqcup}

\define\part{\partial}
\define\emp{\emptyset}

\define\n{\notin}

\define\m{\mapsto}
\define\do{\dots}

\define\sub{\subset}    

\define\T{\times}
\define\ti{\tilde}
\define\nl{\newline}
\redefine\i{^{-1}}

\define\ot{\otimes}
\define\bbq{\bar{\QQ}_l}

\define\tr{\text{\rm tr}}

\redefine\c{\chi}

\redefine\d{\delta}
\define\e{\epsilon}

\define\ph{\phi}

\define\r{\rho}
\define\s{\sigma}
\redefine\t{\tau}
\define\th{\theta}

\redefine\l{\lambda}
\define\z{\zeta}
\define\x{\xi}

\redefine\G{\Gamma}

\redefine\L{\Lambda}
\define\Ph{\Phi}

\define\kk{\bold k}

\define\CC{\bold C}

\define\NN{\bold N}

\define\QQ{\bold Q}

\define\ZZ{\bold Z}
\define\XX{\bold X}

\define\cb{\Cal B}
\define\cc{\Cal C}

\define\ce{\Cal E}

\define\cu{\Cal U}

\define\tE{\ti E}

\define\tS{\ti S}

\head Introduction\endhead
\subhead 0.1\endsubhead
Let $W$ be a Weyl group and let $I$ be the set of simple
reflections of $W$. Let $w$ be the product of the elements of $I$ in
some order (a Coxeter element). Let $\Irr(W)$ be the set of (isomorphism
classes of) irreducible $\CC[W]$-modules. It has been stated by
I.G.Macdonald (around 1970, but unpublished as far as we know) that 

(a) if $E\in\Irr(W)$, then $\tr(w,E)\in\{0,1,-1\}$.
\nl
It is enough to prove (a) when $W$ is irreducible (in the
remainder of this paper we assume that $W$ is irreducible).
If in addition
$W$ is of exceptional type, then (a) can be extracted from the known
character tables. A proof for $W$ of classical type is given in \S1.

In \S2 we describe a connection of (a) with the theory of
two-sided cells \cite{KL} in $W$.

In \S3 we give an Iwahori-Hecke algebra version of (a).

In \S4 we consider the noncrystallographic case.

\subhead 0.2\endsubhead
Let $h$ be the order of $w$ in $W$.
It is known that the centralizer of $w$ in $W$ is a cyclic group of
order $h$.

Assuming that we have found $E_1,E_2,\do,E_h$ distinct in $\Irr(W)$
such that $\tr(w,E_j)=\pm1$ for $j\in[1,h]$, it would follow that
$\tr(w,E)=0$ for any $E\in\Irr(W)$ such that $E\n\{E_1,E_2,\do,E_h\}$.
Indeed, we have $\sum_j\tr(w,E_j)^2=h$; since the centralizer of $w$ in
$W$ is a cyclic group of $h$ we also have
$\sum_{E\in\Irr(W)}\tr(w,E)^2=h$ hence
$\sum_{E\in\Irr(W);E\n\{E_1,\do,E_h\}}\tr(w,E)^2=0$ and our assertion
follows.

Let $r=|I|$ and let $\L^i\in\Irr(W),(i\in[0,r])$ be the $i$-th exterior
power of the reflection representation of $W$.

\subhead 0.3\endsubhead
I thank Eric Sommers for useful comments.

\head 1. Proof of 0.1(a)\endhead
\subhead 1.1\endsubhead
In this subsection $W$ is the symmetric group $S_n$ with the standard
Coxeter group structure. We have $r=n-1$. In this case $w$ is an
$n$-cycle. For $i\in[0,r]$, $\L^i$ correponds to the partition
$1\le1\le1\le\do\le1\le n-i$ of $n$ in the standard parametrization
\cite{L84,4.4} of $\Irr(W)$. We have clearly $\tr(w,\L^i)=(-1)^i$.
By an argument in 0.2 we see that $\tr(w,E)=0$ if
$E\in\Irr(W)-\{\L^0,\do,\L^r\}$. This proves
0.1(a) in our case.

\subhead 1.2\endsubhead
We preserve the setup of 1.1. Assume that $n\ge3$.
Let $w'$ be an $(n-1)$-cycle in $S_n$.

For $i=1,..,n-3$ let ${}'\L^i$ be the irreducible representation
of $S_n$ corresponding to the partition
$1\le1\le1\le\do\le1\le2\le n-i-1$ of $n$
in the standard parametrization of $\Irr(W)$.
We set ${}'\L^0=\L^0$, ${}'\L^{n-2}=\L^{n-1}$.
We have the following result.

(a) $\tr(w',{}'\L^i)=(-1)^i$ for $i\in[0,n-2]$.
Moreover if $E\in\Irr(W)$, $E\n\{{}'\L^0,\do,{}'\L^{n-2}\}$ then
$\tr(w',E)=0$.
\nl
This is proved by induction on $n$, using the results in 1.1 for
$S_{n-1}$ and the known results about restricting an
$E\in\Irr(S_n)$ to $S_{n-1}$.

\subhead 1.3\endsubhead
Let $W_n,n\ge1$ be the group of all permutations of
$\{1,2,\do,n,n',\do,2',1'\}$ which commute with the involution
$i\m i',i'\m i$ ($i\in[1,n]$).
For $i=1,\do,n-1$ let $s_i\in W_n$ be the permutation which
interchanges $i$ with $i+1$ and $i'$ with $(i+1)'$ and leaves the
other elements unchanged.
Let $s_n\in W_n$ be the permutation which interchanges $n$ with $n'$
and leaves the other elements unchanged. Let $I=\{s_1,\do,s_n\}$.
Then $(W_n,I)$ is a Coxeter group of type $B_n$.

In this subsection we assume that $(W,I)=(W_n,I)$.
In this case we have $r=n$ and $w$ is a $2n$-cycle.

Let
$\c:W_n@>>>\{\pm1\}$ be the homomorphism defined by the conditions
$\c(s_i)=1$ for $i=1,\do,n-1$, $\c(s_n)=-1$. A permutation in $W_n$
defines a permutation of the $n$-element set consisting of
$(1,1'),(2,2'),\do,(n,n')$. Thus we have a natural
homomorphism $W_n@>>>S_n$. This induces an injective map
$\Irr(S_n)@>>>\Irr(W_n)$ denoted by $E\m\tE$.
Let $\L^i_n$ ($i=0,\do,n-1$) be the $i$-th exterior power of the
reflection representation of $S_n$. Then, in the parametrization
\cite{L84,4.5} of $\Irr(W_n)$, $\tE_i$ corresponds to the symbol
$$\mat 1&2&3&\do&i&n\\0&1&2&\do&i-1&{}\endmat$$
and $\tE_i\ot\c$  corresponds to the symbol
$$\mat 0&1&2&\do&i-1&i&i+1\\1&2&3&\do&i&n&{}\endmat.$$
From the definition we have (using results in 1.1):

$\tr(w,\tE_i)=(-1)^i$, $\tr(w,\tE_i\ot\c)=(-1)^{i+1}$,
\nl
Since $2n=h$, by an argument in 0.2 we see $\tr(w,E)=0$ if

$E\in\Irr(W)-\{\tE_0,\do,\tE_{n-1},\tE_0\ot\c,\do,\tE_{n-1}\ot\c\}$.
\nl
This proves 0.1(a) in our case.

\subhead 1.4\endsubhead
In this subsection we assume that $n\ge4$
and that $W=W'_n:=\ker(\c)\sub W_n$ (notation of 1.3) with $I$
being the subset of $W'_n$ consisting of $s_i (i=1,\do,n-1)$
and of $s_ns_{n-1}s_n$. (This $(W,I)$ is a Coxeter group of type
$D_n$.)
In this case we have $r=n$.

We can identify $W_{n-1}\T W_1$ with the subgroup of $W_n$
consisting of all permutations which map
$\{1,2,\do,n-1,(n-1)',\do,2',1'\}$ to itself and hence also map
$\{n,n'\}$ to itself. We can assume that $w$ is the product of a
$(2n-2)$-cycle in $W_{n-1}$ with a $2$-cycle in $W_1$.
As in 1.3 we have an obvious surjective
homomorphism  $W_{n-1}\T W_1@>>>S_{n-1}\T S_1=S_{n-1}$.
Via this homomorphism we can regard $\L^i_{n-1}$ (with $i\in[0,n-2]$), 
see 1.3, as a $W_{n-1}\T W_1$ module; inducing this from
$W_{n-1}\T W_1$ we obtain an irreducible
$W_n$-module of $W_n$ whose restriction to $W'_n$
is an irreducible $W'_n$-module $M_i$ which in the parametrization
\cite{L84,4.6} of $\Irr(W'_n)$ corresponds to the symbol
$$\mat1&2&3&\do&i-1&i&n-1\\0&1&2&\do&i-2&i-1&i+1\endmat.
$$
From the definition we have (using results in 1.1):

$\tr(w,M_i)=(-1)^i$.

For $i\in[0,n-2]$ the irreducible $S_n$-module ${}'\L_i$ (see 1.2)
can be regarded as a $W_n$-module; restricting this to $W'_n$ gives
a $W'_n$-module ${}'M_i$ which in the parametrization \cite{L84,4.6}
of $\Irr(W'_n)$ corresponds to the symbol
$$\mat1&2&3&\do&i-1&i+1&n-1\\0&1&2&\do&i-2&i-1& i\endmat$$
if $i\in[1,n-3]$, to the symbol
$$\mat1&2&3&\do&n-1&n\\0&1&2&\do&n-2&n\endmat$$
if $i=n-2$ and to the symbol
$$\mat n\\0\endmat$$
if $i=0$.
From the definition we have (using results in 1.2):

$\tr(w,{}'M_i)=(-1)^i$.

Since $2n-2=h$, by an argument in 0.2 we see that $\tr(w,E)=0$ if
$$E\in\Irr(W)-\{M_0,\do,M_{n-2},{}'M_0,\do,{}'M_{n-2}\}.$$
This proves 0.1(a) in our case.

\head 2. Relation with two-sided cells\endhead
\subhead 2.1\endsubhead
In this subsection $W$ is of type $B_n$. Let $i\in[0,n]$.
In the parametrization \cite{L84,4.5} of $\Irr(W)$, $\L^i$
corresponds to the symbol
$$\mat0&1&2&\do&i-1&n\\1&2&3&\do&i&{}\endmat.$$

\subhead 2.2 \endsubhead
In this subsection $W$ is of type $D_n,n\ge4$. 
In the parametrization \cite{L84,4.6} of $\Irr(W)$, $\L^j$
(with $j\in[1,n]$) corresponds to the symbol
$$\mat1&2&3&\do&j-1&j\\0&1&2&\do&j-2&n-1\endmat$$
and $\L^0$ corresponds to the symbol  
$$\mat n\\0\endmat.$$

\subhead 2.3 \endsubhead
Let $\cc$ be the set of two-sided cells of $W$. 

It is known that we have a well defined surjective map
$u:\Irr(W)@>>>\cc$. For $i\in[0,r]$ let $c_i=u(\L^i)\in\cc$. 

For $c\in\cc$ we denote by $\G_c$
the finite group associated to $c$ in \cite{L84}.

For $i\in[0,r]$, $\G_{c_i}$ is one of the symmetric groups
$S_1,S_2,S_3,S_4,S_5$. More precisely, we describe the sequence
$(\G_{c_0},\G_{c_1},\do,\G_{c_r})$ for the various $W$. For $W$
simplylaced this sequence is as follows.

Type $A$: $S_1,S_1,\do,S_1$;

Type $D$: $S_1,S_1,S_2,S_2,\do, S_2,S_1,S_1$;

Type $E_6$: $S_1,S_1,S_2,S_3,S_2,S_1,S_1$;

Type $E_7$: $S_1,S_1,S_2,S_3,S_3,S_2,S_1,S_1$;

Type $E_8$: $S_1,S_1,S_2,S_3,S_5,S_3,S_2,S_1,S_1$.

For $W$ non-simplylaced this sequence is as follows.

Type $B=C$: $S_1,S_2,S_2,\do, S_2,S_1$;

Type $F_4$: $S_1,S_2,S_4,S_2,S_1$;

Type $G_2$: $S_1,S_3,S_1$.

When $W$ is of type $A$ this is obvious. When $W$ is of type $B$
or $D$ this follows from 2.1, 2.2. When $W$ is of exceptional type
this follows from the tables in the Appendix to \cite{L84}.

\subhead 2.4\endsubhead
Let $\cc_{ex}$ be the set of exceptional two-sided cells of $W$
in the sense of \cite{L17,1.3}. This set is empty if $W$ is not of type
$E_7,E_8$. If $W$ is of type $E_7$, $\cc_{ex}$ consists of a single
cell denoted by $c^{7/2}$. If $W$ is of type $E_8$, $\cc_{ex}$ consists
of two cells; one of them contains the cell $c^{7/2}$ of $E_7$ and is
denoted by $c^{9/2}$; the other one is denoted by $c^{7/2}$.

\subhead 2.5\endsubhead
Let $\CC[\Irr(W)]$ be the free $\CC$-vector spacee with basis
$\{E;E\in\Irr(W)\}$.
We set $\ce(W)=\sum_{E\in\Irr(W)}\tr(w,E)E\in\CC[\Irr(W)]$.
We have $\ce(W)=\sum_{c\in\cc}\ce_c(W)$ where
$$\ce_c(W)=\sum_{E\in\Irr(W);u(E)\in c}\tr(w,E)E.\tag a$$
For $c\in\cc$ let
$\XX_c=\{E\in\Irr(W);u(E)\in c,\tr(w,E)\ne0\}$
and let $n_c$ be the number of elements of $\XX_c$
that is the number of nonzero terms in the sum (a).

We have the following result.

\proclaim{Proposition 2.6}(a) Assume that $c=c_i$ for some $i\in[0,r]$
($i$ is necessarily unique). If $\G_{c_i}=S_1$ then $n_{c_i}=1$. If
$\G_{c_i}=S_2$ then $n_{c_i}=2$. If $\G_{c_i}=S_3$ then $n_{c_i}=4$. If
$\G_{c_i}=S_4$ then $n_{c_i}=6$. If $\G_{c_i}=S_5$ then $n_{c_i}=10$.

(b) Assume that $c\in\cc_{ex}$. Then $n_c=2$.

(c) Assume that $c\in\cc$ is not of the form $c_i$ for some $i$ and
that $c\n\cc_{ex}$. Then $n_c=0$. In particular, we have $\ce_c(W)=0$.
\endproclaim
In the case where $W$ is of type $A$ this follows from 1.1.
In the case where $W$ is of type $B$ this follows from 1.3 and 2.1.
In the case where $W$ is of type $D$ this follows from 1.4 and 2.2.
In the case where $W$ is of exceptional type this follows by
examining the existing tables.

\subhead 2.7\endsubhead
Assume that $W$ is of type $B_4$.
We describe the elements $\ce_{c_i}(W)$ for $i=0,1,2,3,4$. (We identify
objects in $\Irr(W)$ with the corresponding symbols.)

$\mat4\\-\endmat$,

$-\mat1&4\\0\endmat-\mat0&1\\4\endmat$,

$\mat1&2&4\\0&1\endmat+\mat0&1&2\\1&4\endmat$,

$-\mat1&2&3&4\\0&1&2\endmat-\mat0&1&2&3\\1&2&4\endmat$,

$\mat0&1&2&3&4\\1&2&3&4\endmat$.

\subhead 2.8\endsubhead
Assume that $W$ is of type $D_5$.
We describe the elements $\ce_{c_i}(W)$ for $i=0,1,2,3,4,5$. (We
identify objects in $\Irr(W)$ with the corresponding symbols.)

$\mat5\\0\endmat$,

$-\mat4\\1\endmat$,

$\mat1&4\\0&2\endmat-\mat2&4\\0&1\endmat$,

$-\mat1&2&4\\0&1&3\endmat+\mat1&3&4\\0&1&2\endmat$,

$\mat1&2&3&4\\0&1&2&4\endmat$,

$-\mat1&2&3&4&5\\0&1&2&3&4\endmat$.

\subhead 2.9\endsubhead
If $W$ is an exceptional Weyl group and $E\in\Irr(W)$ we write $E=d_b$
when $\dim(E)=d$ and $b$ is the smallest integer such that $E$ appears
in the $b$-th symmetric power of the reflection representation. (When
there are two such $E$ we denote them by $d_b,d'_b$.)

\subhead 2.10\endsubhead
Assume that $W$ is of type $E_6$.
We describe the elements $\ce_{c_i}(W)$ for $i=0,1,2,3,4,5,6$.

$1_0$,

$-6_1$,

$30_3-15_4$,

$60_8-90_8-10_9+20_{10}$,

$30_{15}-15_{16}$,

$-6_{25}$,

$1_{36}$.

\subhead 2.11\endsubhead
Assume that $W$ is of type $E_7$.
We describe the elements $\ce_{c_0}(W)$, $\ce_{c_1}(W)$, $\ce_{c_2}(W)$,
$\ce_{c_3}(W)$, $\ce_{c^{7/2}}(W)$, $\ce_{c_4}(W)$, $\ce_{c_5}(W)$,
$\ce_{c_6}(W)$, $\ce_{c_7}(W)$.

$1_0$,

$-7_1$,

$56_3-35_4$,

$280_8-280_9-70_9+35_{13}$,

$512_{11}-512_{12}$,

$-280_{17}+280_{18}+70_{18}-35_{22}$,

$-56_{30}+35_{31}$,

$7_{46}$,

$-1_{63}$.

\subhead 2.12\endsubhead
Assume that $W$ is of type $E_8$.
We describe the elements $\ce_{c_0}(W)$, $\ce_{c_1}(W)$, $\ce_{c_2}(W)$,
$\ce_{c_3}(W)$, $\ce_{c^{7/2}}(W)$, $\ce_{c_4}(W)$, $\ce_{c^{9/2}}(W)$,
$\ce_{c_5}(W)$, $\ce_{c_6}(W)$, $\ce_{c_7}(W)$, $\ce_{c_8}(W)$.

$1_0$,

$-8_1$,

$112_3-84_4$,

$1344_8-1008_9-448_9+56_{10}$,

$4096_{11}-4096_{12}$,

$$\align&4480_{16}-7168_{17}-5670_{18}+2016_{19}+1134_{20}\\&
-420_{20}+1680_{22}+4536_{23}-448_{25}-70_{32},\endalign$$

$-4096_{26}+4096_{27}$,

$1344_{38}-1008_{39}-448_{39}+56_{40},$

$112_{63}-84_{64},$

$-8_{91}$,

$1_{120}$,

\subhead 2.13\endsubhead
Assume that $W$ is of type $F_4$.
We describe the elements $\ce_{c_i}(W)$, $i=0,1,2,3,4$.

$1_0$,

$-2_4-2'_4$,

$12_4-6_6-6'_6+4_8+1_{12}+1'_{12}$,

$-2_{16}-2'_{16}$,

$1_{24}$.

\subhead 2.14\endsubhead
Assume that $W$ is of type $G_2$.
We describe the elements $\ce_{c_i}(W)$, $i=0,1,2$.

$1_0$,

$-2_1+2_2+1_3+1'_3$,

$1_6$.

\subhead 2.15\endsubhead
The pattern described in 2.7-2.14 is similar to that
of unipotent representations of a finite Chevalley group
associated to a Coxeter orbit
in \cite{L76}; in fact the first can be deduced from the second
using the nonabelian Fourier transform \cite{L84,4.14}.

\subhead 2.16\endsubhead
Let $E\in\Irr(W)$. Let $f_E\in\ZZ_{>0}$ be as in \cite{L84, 4.1.1}.
For $i\in[0,r]$ let
$$p_i=\sum_{E\in u\i(c_i)}\tr(w,E)1/f_E.$$
It turns out that $p_i$ is always an integer (we don't have an
explanation of this fact.)

If $\G(c_i)=S_1$ then $p_i=(-1)^i$.

If $\G(c_i)=S_2$ and $W$ is simply laced then $p_i=0$.
(We use that $1/2-1/2=1$.)

If $\G(c_i)=S_2$ and $W$ is not simply laced then $p_i=(-1)^i$. (We
use that $1/2+1/2=1$.)

If $\G(c_i)=S_3$ and $W$ is simply laced then $p_i=0$. (We use that
$1/2-1/3-1/3+1/6=0$.)

If $\G(c_i)=S_3$ and $W$ is not simply laced then $p_i=1$. (We use that
$-1/6+1/2+1/3+1/3=1$.)

If $\G(c_i)$ is $S_4$ or $S_5$ then $p_i=0$. (We use that

$1/24-1/12-1/3+1/8+1/8+1/8=0$,

$1/120-1/12-1/30+1/6+1/6-1/5+1/20+1/24-1/12-1/30=0$.)

\subhead 2.17\endsubhead
For $i\in[0,r]$ let $\r_i$ be the unique special representation in
$u\i(c_i)$.

If $\G(c_i)=S_1$ then $\r_i=\L^i\in\XX_{c_i}$.

If $\G(c_i)=S_2$ and $W$ is simply laced then $\r_i\in\XX_{c_i}$,
$\L^i\n\XX_{c_i}$.

If $\G(c_i)=S_2$ and $W$ is not simply laced then $\r_i\n\XX_{c_i}$,
$\L^i\n\XX_{c_i}$.

If $\G(c_i)=S_3$ and $W$ is simply laced then $\r_i\n \XX_{c_i}$,
$\L^i\in\XX_{c_i}$.

If $\G(c_i)=S_3$ and $W$ is of not simply laced then $\r_i\in\XX_{c_i}$,
$\L^i\in\XX_{c_i}$.

If $\G(c_i)$ is $S_4$ or $S_5$ then $\r_i\in\XX_{c_i}$,
$\L^i\in\XX_{c_i}$.

\head 3. An Iwahori-Hecke algebra analogue of 0.1(a)\endhead
\subhead 3.1\endsubhead
Let $l:W@>>>\NN$ be the standard length function. Let $w_0$ be the
longest element of $W$ and let $\nu=l(w_0)$.
Let $v$ be an indeterminate. Let $K=\CC(v)$. Let $H$ be the $K$-vector
space with basis $\{T_w;w\in W\}$. It is known that $H$ has a well
defined associative $K$-algebra multiplication with unit $T_1=1$
in which

$T_{s_i}^2=(v^2-1)T_{s_i}+v^2$ for any $i\in I$,

$T_yT_{y'}=T_{yy'}$ for any $y,y'$ in $W$ such that
$l(yy')=l(y)+l(y')$.

$H$ is the Iwahori-Hecke algebra associated to $W$.
It is known that to each $E\in\Irr(W)$ one can associate
canonically an $H$-module $E(v)$ such that for any $y\in W$ we have
$\tr(T_y,E(v))\in\CC[v]$ and $\tr(T_y,E(v))_{v=1}=\tr(y,E)$.
We have the following version of 0.1(a).

(a) If $E\in\Irr(W)$, then $\tr(T_w,E(v))=\tr(w,E)v^{m_E}$
for some $m_E\in\NN$. In particular $\tr(T_w,E(v))$ is either $0$
or $\pm v^{m_E}$.

\subhead 3.2\endsubhead
It is known that $\tr(T_w,E(v))$ (for $E\in\Irr(W)$) is independent of
the choice
of $w$ in 0.1. It is also known that $w$ can be chosen so that
the following holds:

(a) If $W$ is of type $\ne A$ or of type $A$ with $r$ odd (so that
$h$ is even) then $w^{h/2}=w_0$, $l(w)h/2=l(w_0)$.

(b) If $W$ is of type $A$ with
$r$ even (so that $h=r+1$ is odd) then for some $a,b$ in $W$
we have $w=ab$, $l(a)=l(b)=r/2$, $w_0=abab\do aba=baba\do bab$
($r+1$ factors) and $l(w_0)=(r+1)r/2$.

Thus we may assume that (a) or (b) holds.
If (a) holds we have $T_w^{h/2}=T_{w_0}$. Hence $T_w^h=T_{w_0}^2$.
If (b) holds we have
$T_{w_0}=T_aT_bT_a\do T_aT_bT_a=T_bT_a\do T_bT_aT_b$
(both products involve $r+1$ factors). Hence

$T_{w_0}^2=(T_aT_bT_a\do T_aT_bT_a)(T_bT_a\do T_bT_aT_b)=
(T_aT_b)^{r+1}=T_w^h$.

Thus in both cases we have $T_{w_0}^2=T_w^h$.

Let $E\in\Irr(W)$. Let $a_E\in\NN$ be as in \cite{L84, 4.1.1}.
Let $E^!=E\ot\L^r$.

By \cite{L84, 5.12.2} there is a $K$-linear involution
$\s:E(v)@>>>E(v)$ such that $\s T_{w_0}=T_{w_0}\s=v^{\nu-a_E+a_{E^!}}$
on $E(v)$. We then have
$T_{w_0}^2=v^{2(\nu-a_E+a_{E^!})}$ on $E(v)$.
Since $T_w^h=T_{w_0}^2$, it follows that any eigenvalue of
$T_w:E(v)@>>>E(v)$ (in an algebraic closure $\bar K$ of $K$)
is of the form $\x Z$ where $\x\in\CC^*,\x^h=1$ and $Z$ is a fixed
element of $\bar K$ such that $Z^h=v^{2(\nu-a_E+a_{E^!})}$.
It follows that $\tr(T_w,E(v))=z_E Z$ for some $z_E\in\CC$.
Assume now that $z_E\ne0$. Since $\tr(T_w,E(v))\in\CC[v]$,
it follows that $2(\nu-a_E+a_{E^!})=hm_E$ for some $m_E\in\NN$
so that $Z=\z v^{m_E}$ for some $\z\in\CC^*,\z^h=1$ and
$\tr(T_w,E(v))=z'_E v^{m_E}$ for some $z'_E\in\CC^*$. Specializing
$v=1$ we obtain $\tr(w,E)=z'_E$ and by 0.1(a) we have $z'_E=\pm1$.
This proves 3.1(a). (A similar argument appears in \cite{L84,p/320}.)

We see also that in 3.1(a) we have

(c) $m_E=2(\nu-a_E+a_{E^!})/h$ whenever $\tr(w,E)\ne0$.
\nl
From (c) we can calculate $m_E$ in each case; we obtain the
following result.

\proclaim{Proposition 3.3}If $E\in\XX_{c_i}$ (see 2.5) with $i\in[0,r]$,
or if $E\in\XX_{c^i}$ (see 2.4, 2.5) with $i\in\{7/2,9/2\}$, then
$m_E=2r-2i$.
\endproclaim

\head 4. The noncrystallographic case\endhead
\subhead 4.1\endsubhead
In this section $(W,I)$ is an (irreducible) noncrystallographic finite
Coxeter group. Let $w,\Irr(W),r,\L^i (i\in[0,r])$
be objects analogous to those in 0.1, 0.2.
In this case the analogue of 0.1(a) does not hold.

Let $\cc$, $u:\Irr(W)@>>>\cc$, $c_i$ be objects analogous
to those in 2.3.

Let $\cc_{ex}$ be the set of exceptional two-sided cells of $W$
in the sense of \cite{L17,1.3}. This set is empty if $W$ is not of type
$H_3,H_4$. If $W$ is of type $H_3$, $\cc_{ex}$ consists of a single
cell denoted by $c^{3/2}$. If $W$ is of type $H_4$, $\cc_{ex}$
consists of two cells; one of them contains the cell $c^{3/2}$ of $H_3$
and is denoted by $c^{5/2}$; the other one is denoted by $c^{3/2}$.

Let $\ce(W)\in\CC[\Irr(W)], \ce_c(W)\in\CC[\Irr(W)]$
($c\in\cc$) be the objects analogous to those in 2.5. We have the
following analogue of 2.6.

(a) For $c\in\cc$ we have $\ce_c(W)\ne0$ if and only if $c=c_i$ for
some $i\in[0,r]$ or $c=c^i$ for some $i\in\{3/2,5/2\}$.

The results in 4.2-4.5 are extracted from the known character tables
of $W$. For $E\in\Irr(W)$ we use a notation analogous to that in 2.9.

\subhead 4.2\endsubhead
Assume that $W$ is of type $H_3$. We set $\z=e^{2\pi\sqrt{-1}/5}$.
We describe the elements $\ce_{c_0}(W)$, $\ce_{c_1}(W)$,
$\ce_{c^{3/2}}(W)$, $\ce_{c_2}(W)$, $\ce_{c_3}(W)$.

$1_0$,

$(\z+\z^4)3_1+(\z^2+\z^3)3_3$,

$4_3-4_4$,

$(\z+\z^4)3_6+(\z^2+\z^3)3_8,$

$1_{15}$.

\subhead 4.3\endsubhead
Assume that $W$ is of type $H_4$. We set $\z=e^{2\pi\sqrt{-1}/5}$.
We describe the elements $\ce_{c_0}(W)$, $\ce_{c_1}(W)$,
$\ce_{c^{3/2}}(W)$, $\ce_{c_2}(W)$, $\ce_{c^{5/2}}(W)$, $\ce_{c_3}(W)$,
$\ce_{c_4}(W)$.

$1_0$,

$(\z+\z^4)4_1+(\z^2+\z^3)4_7$,

$16_3-16_6$,

$$\align&-(\z^2+\z^3)24_6+24_7+40_8-48_9+(\z+\z^4)30_{10}
+\\&(\z^2+\z^3)30'_{10}+(\z^2+\z^3)16_{11}+24_{11}-(\z+\z^4)6_{12}\\&
-(\z+\z^4)24_{12}-10_{12}-8_{13}+(\z+\z^4)16_{13}-(\z^2+\z^3)6_{20},
\endalign$$

$-16_{18}+16_{21}$,

$(\z+\z^4)4_{31}+(\z^2+\z^3)4_{37}$,

$1_{60}$.

\subhead 4.4\endsubhead
Assume that $W$ is a dihedral group of order $2(2p+1)$.

We set $\z=e^{2\pi\sqrt{-1}/(2p+1)}$.
We describe the elements $\ce_{c_i}(W)$, $i=0,1,2$.

$1_0$,

$(\z+\z^{2p})2_1+(\z^2+\z^{2p-1})2_2+\do+(\z^p+\z^{p+1})2_p$,

$1_{2p+1}$.

\subhead 4.5\endsubhead
Assume that $W$ is a dihedral group of order $2(2p)$.
We set $\z=e^{\pi\sqrt{-1}/p}$.
We describe the elements $\ce_{c_i}(W)$, $i=0,1,2$.

$1_0$,

$(\z+\z^{2p-1})2_1+(\z^2+\z^{2p-2})2_2
+\do+(\z^{p-1}+\z^{p+1})2_{p-1}-1_p-1'_p$,

$1_{2p}$.

Note that if $p$ is even then the coefficient of $2_{p/2}$ above is zero
so that $2_{p/2}$ does not actually appear in the sum.

\subhead 4.6\endsubhead
The first assertion of 3.1(a) continues to hold in our case, with
the same proof.

\enddocument